\documentclass[a4paper,11pt]{amsart}
\usepackage[utf8]{inputenc}
\usepackage{amsmath, amsthm, amscd, amsfonts, amssymb, graphicx, color}
\usepackage[bookmarksnumbered, colorlinks, plainpages]{hyperref}
\usepackage{wrapfig}
\usepackage{subfig}
\usepackage{subfig}
\usepackage{cite}
\usepackage[framemethod=default]{mdframed}
\usepackage{showexpl}
\usepackage[algoruled]{algorithm2e}

\begin{document}
	
	 \newcommand{\be}{\begin{equation}}
	 \newcommand{\ee}{\end{equation}}
	 \newcommand{\bt}{\beta}
	 \newcommand{\al}{\alpha}
	 \newcommand{\laa}{\lambda_\alpha}
	 \newcommand{\lab}{\lambda_\beta}
	 \newcommand{\no}{|\Omega|}
	 \newcommand{\nd}{|D|}
	 \newcommand{\Om}{\Omega}
	 \newcommand{\h}{H^1_0(\Omega)}
	 \newcommand{\lt}{L^2(\Omega)}
	 \newcommand{\la}{\lambda}
	 \newcommand{\ro}{\varrho}
	 \newcommand{\cd}{\chi_{D}}
	 \newcommand{\cdc}{\chi_{D^c}}
	 \newtheorem{thm}{Theorem}[section]
	 \newtheorem{cor}[thm]{Corollary}
	 \newtheorem{lem}[thm]{Lemma}
	 \newtheorem{prop}[thm]{Proposition}
	 \theoremstyle{definition}
	 \newtheorem{defn}{Definition}[section]
	 \newtheorem{exam}{Example}[section]
	 \theoremstyle{remark}
	 \newtheorem{rem}{Remark}[section]
	 \numberwithin{equation}{section}
	 \renewcommand{\theequation}{\thesection.\arabic{equation}}
	 \numberwithin{equation}{section}
	\newtheorem{alg}{Algorithm}

\newcommand{\R}{\mathbb{R}}
\renewcommand{\div}{\operatorname{div}}
\newcommand{\st}{\;:\;}
\newcommand{\dx}{\,\mathrm{d}x}
\renewcommand{\d}{\mathrm{d}}
\newcommand{\norm}[1]{\left\|#1\right\|}
\newcommand{\mean}{\operatorname{mean}}

	\title[ Numerical approximations  of chromatographic models]{Numerical approximations  of chromatographic models}
	\author[Farid Bozorgnia,  Sonia Seyed Allaei   ]{ Farid Bozorgnia,  Sonia Seyed Allaei    }
	 \address{CAMGSD, Instituto Superior T\'{e}cnico, University of Lisbon, Av. Rovisco Pais, 1049-001 Lisbon, Portugal.}
 \email{farid.bozorgnia@tecnico.ulisboa.pt,\,
   Sonia.seyedallaei@tecnico.ulisboa.pt}

		 \date{\today}
	
	 \thanks{*The corresponding author, F. Bozorgnia was  supported by the Portuguese National Science Foundation through FCT fellowships SFRH/BPD/33962/2009 and  project UTAP-EXPL/MAT/0017/2017.}
	
	 \begin{abstract}
	 	A numerical scheme based on modified method of characteristics with adjusted advection (MMOCAA) is proposed to approximate the solution of   the system  liquid chromatography  with   multi components case.  For the case of one component,  the method preserves the mass. Various examples and computational tests  numerically verify  the accuracy and efficiency  of the approach.
	 \end{abstract}

	 \maketitle

	 \textbf{Keywords}:  Advection-Diffusion, Coupled system,  Langmuir adsorption model,  Liquid Chromatography, Numerical approximation.  \\
	 2010 MSC:

\section{Introduction and problem setting}

Chromatography is a technical   process to separate mixed chemical components  with   a wide range of chemical industrial
applications  such as in  pharmaceutical, food ingredients, etc. Here,  we briefly explain the separation of components
by liquid chromatography. In column chromatography a mixed sample is injected into a fluid stream  which is called mobile phase. Then the fluid  is pumped through a pipe which we refer  as chromatographic column. The column is  filled with very small  porous beads called stationary phase.
Different  components in fluid adsorbs and/or desorbs at different rates on the stationary phase so they move through the  column  at different speeds  and exit the column  at different times;  elution,   see \cite{GG, Js}.

The transport of solutes in heterogeneous porous media is described by mass balance equation.  The transport is influenced by  the convection, diffusion,   dispersion and also reaction/adsorption  between  solute and the porous environment. The model consists of system of convection-diffusion-reaction partial
differential equations with dominating the convective terms coupled via differential or algebraic
equations. To see different models  and numerical approach,  we refer to \cite{Js}.

\section{Preliminaries and Problem setting}
 In one-dimension, the transport is given by the following coupled equations:

	\begin{equation}\label{s1}
\left\{
\begin{array}{ll}
\frac{\partial u_{i}}{\partial  t } +  F\,    \frac{\partial q_{i}}{\partial  t }  +
  v  \frac{\partial u_{i}}{\partial  x }  = D  \, \frac{\partial^{2}  u_{i}}{\partial  x^2}  	 &	 \text{in }(0,L)\times(0, T),\\\\
u_{i}( x=0, t) = g_{i}(t)  	&	\text{on } (0,T),\\\\
u_{i} (x, t=0)=  u_{i,0}(x) & \text{in } (0, L),\\\\
\frac{\partial u_{i}(x=L,t) }{\partial  x }=0   & \text{for  } t\ge 0,\\\\
i=1,2, \cdots m,
\end{array}
\right.
\end{equation}
  where,
\begin{itemize}
\item $L:$ the column length,
\item   $t:$ time,
\item   $u_{i}:$  concentration  of      the $i^{\rm{th}}$   component in the mobile phase,
\item    $q_{i}:$   concentration  of  the $i^{\rm{th}}$ component     in the stationary phase,
  \item  $u_{i,0}:$      initial condition,
  \item    $g_{i}:$    boundary condition (injection profile),
\item  $F:$      stationary/mobile phase ratio,
  \item  $v:$    mobile phase velocity,
       \item $D:$   diffusion parameter,
\item   $m:$    number of mixture components in the sample.

   \end{itemize}
   The Neumann  boundary condition   $\frac{\partial u_{i}(x=L,t) }{\partial  x }=0$  persuade continuity of the outlet concentration profile to the connecting tube receiving the fluid after leaving the column.
The dispersion coefficient $D$  is given by
\[
D= \frac{L\, v}{ 2 N_t},
\]
where $N_t$ is the number of theoretical plates.  The term $F$ is given by
\[
F=\frac{1 -  \epsilon}{\epsilon},
\]
 which indicates  the phase ratio based on the porosity $\epsilon$. Also  $q_i$ is called adsorption isotherm and  we assume that $q_i=q_i(u_1, \cdots, u_m)$.  In  Langmuir model this term is given by
\begin{equation}\label{Q}
q_{i} = \frac{a_{i}u_{i}}{1+ \sum\limits_{j=1}^{m} b_{j} u_{j}}, \quad i=1,2, \cdots m,
\end{equation}
where $ a_i, b_i >0$.


Let  assume that the  mass  of components at the initial time in column is zero; $u_{i,0}(x)=0,\,\,\,\,0<x< L$.
 We consider  rectangular injection profiles so boundary condition at the inlet point is:
\[
u_{i}(x=0, t) = g_i(t) = u_{i,\textrm{inj}}(t),
\]
with
	\begin{equation}\label{ss1}
u_{i,\textrm{inj}} (t)=
\left\{
\begin{array}{ll}
u_{i,f}(t)    &  0 < t \le t_{\textrm{inj}},\\
0          &   t>  t_{\text{inj}},
\end{array}
\right.
\end{equation}
where  $u_{i,f}(t)$ is the inlet feed concentration and  $t_{\text{inj}}$ is the injection time. One can consider Danckwerts-type  boundary conditions  at the column  inlet which is given by,
\[
u_{i}( x=0, t)=u_{i,\textrm{inj}}(t)+ \frac{D}{v}\, \frac{\partial u_i}{\partial x} \quad i=1,\cdots ,m,
\]
where for  $N_t > 100$, e.g. Seidel-Morgenstern \cite{S1}, it reduces again to
\[
 u_{i}( x=0, t)=u_{i,\text{inj}}(t).
 \]
%

It is well known that in the convection dominate problems,   discontinuity propagates in time even with the smooth initial and boundary data. Furthermore, the nonlinearity  and coupling in term $q_{i}$  in (\ref{s1})  brings more challenges  to the numerical solution of  this type of nonlinear coupled convection-diffusion system.

Standard finite difference, finite volume, and finite element methods are not stable and  the numerical   approximations    exhibit non-physical oscillations and/or generates  artificial numerical diffusion, which smear out sharp fronts of the solution \cite{GSK,Hb,RUE}.

In the  case of scalar equation, one   approach   to eliminate the  nonphysical oscillation which occurs on standard finite element or  finite difference approach, is based on   characteristic method. The sketch of idea is  splitting the equation  into two sub-steps, the convection step, which is solved explicitly by high order schemes (Lax  Wendroff for instance), and the diffusion step, which is solved implicitly by central difference, see \cite{B1, D1}.

For the    system    (\ref{s1}),  different approaches have been discussed.
 In \cite{J1}  high resolution semi-discrete flux-limiting
finite volume scheme is proposed   which is capable to defeat   numerical oscillations and  preserves the
positivity of numerical solution.  The authors   validate their scheme
against other flux-limiting schemes available in the literature.     To see about   discontinuous Galerkin  approximation  for system  (\ref{s1})  we refer to   \cite{J2, M1, M2}. Recently in \cite{Q2} a  transport model is used to describe gradient elution in liquid chromatography. Furthermore, the  authors implement  Laplace transform to obtain the   analytical  solution of model.

	
In \cite{B2} the  existence of the unique weak solution has shown for  the case that $q= \nabla \phi$ for  some   $\phi  :\mathbb{R}^m  \rightarrow \mathbb{R}$, i.e  the vector field
  $ q:  \mathbb{R}^m  \rightarrow \mathbb{R}^m$ in (\ref{s1}) can be expressed as a gradient of some non-negative  $C^1$-convex function $\phi$. The proof is based  on Rothe's method along with solving a convex minimization problem at each time step which gives a numerical method to solve the coupled system.

We propose the modified method of characteristics with adjusted advection (MMOCAA) to solve the system of equation (\ref{s1}). This method was proposed by Douglas et al. to solve
advection dominate transport PDEs \cite{D1}. The MMOCAA corrects  the mass error occurs in the modified method of characteristic (MMOC)
by perturbing the foot of the characteristics vaguely \cite{D2,Ew}.  Our method  is straight forward  to implement and robust comparing the other methods mentioned above. Error analysis for presented scheme  is beyond our aim in the current work.	

The paper is organized as follows. Section 2 deals with introducing problem and previous works. In Section 3 we present our numerical scheme  for coupled system and  for scalar equation in ideal case.  We finally represent    various examples and computational tests.

\section{The numerical scheme}
For the sake of simplicity, let's assume that  the number of components  is two ($m=2$) however, it can simply extended for $m>2$.
\begin{equation}\label{s3}
\left\{
\begin{array}{ll}
\frac{\partial u_{1}}{\partial  t } +    F \frac{\partial q_{1}}{\partial  t }  +   v\,   \frac{\partial u_{1}}{\partial  x }  = D \, \frac{\partial^{2} u_{1}}{\partial  x^2 }    	 &	 \text{in }(0,L)\times(0, T),\\\\
\frac{\partial u_{2}}{\partial  t } +    F \frac{\partial q_{2}}{\partial  t }  +   v\,    \frac{\partial u_{2}}{\partial  x }  = D \, \frac{\partial^{2} u_{2}}{\partial  x^2 }   	 &	 \text{in }(0,L)\times(0, T),\\\\
u_{k}(x=0, t) = g_{i}(t),  \quad  k=1,2    	&	\text{on } (0,T),\\\\
u_{k} (x, t=0)= 0, \,   \quad   \quad  k=1,2  & \text{in } (0, L),\\\\
\frac{\partial u_{1}(x=L,t) }{\partial  x }=0,  \frac{\partial u_{2}(x=L,t) }{\partial  x }=0    & \text{for  } t\ge 0.
\end{array}
\right.
\end{equation}

We start semi-discritization in time for system (\ref{s3}).
For positive integer number $N$,    the time interval $[0, T]$  is divide  to $N$ sub interval as
\begin{equation}
[0,T]=[t^0,t^1]\cup\dots  \cup[t^{N-1},t^N],
\end{equation}
where  $t^{n}=n  \triangle t,n=1,\dots,N$ and $ \triangle t= \frac{T}{N}$.

Let $u_{k}^{n}(x) := u_{k}(t^{n}, x)$.
If we start form the point $(x_i,t^{n+1})$ and move back in direction  of characteristic line, then we hit  the time level $n$. The intersection point is called  $(\tilde x_i, t^{n})$.

By method of characteristic we have
\[
 u_{k} (x_{i},  t^{n+1} )= \tilde u_{k}^{n}(\tilde x_i):=  u_{k}^{n} (x_{i} -  v \, \triangle t),\quad \text{for}  \, k=1,2.
\]
As $\tilde x_i$ may not be a grid point, $\tilde u_{i}^{n}(\tilde x_i)$ is an interpolated value.
For  $v\,   \frac{\triangle t}{ \triangle x}  < 1,$  the foot of backward characteristic  $\tilde x_i$  intersects $ t = t^n$  inside the interval   $(x_{i-1} , x_{i+1}).$   We can use quadratic interpolation between $  u_{i-1}^n, u_{i}^n, u_{i+1}^n$ which  leads to the Lax-Wendrof scheme in the scalar case.

By using the chain rule,   we have
\[
\frac{\partial q_{1}}{\partial t}= \frac{\partial q_{1}}{\partial u_{1}} \, \frac{\partial u_{1} }{ \partial t}   + \frac{ \partial q_{1}}{\partial u_{2}} \, \frac{\partial u_{2} }{ \partial t}.
\]
	  We use the notation $\mathbf{u}(\cdot, t^{n}) = \mathbf{u}^{n}(\cdot).$ To update the values of $\mathbf{u}=(u_1, u_2) $ at the point $(x_{i} , t^{n+1}),$ we follow  backward  in the direction of the characteristic line.
The  semi-discretization of (\ref{s3}) reads as follows
\begin{equation}\label{s4}
\frac{u_{1}^{n+1} - \widetilde{u_{1}}^{n}}{\triangle t} + F \, \frac{\partial q_{1}^{n}}{\partial u_{1}}\,\frac{u_{1}^{n+1} -  u_{1}^{n}}{\triangle t }  + F \, \frac{\partial q_{1}^{n}}{\partial u_{2}}\,\frac{u_{2}^{n+1} -  u_{2}^{n}}{\triangle t}    = D (u_{1}^{n+1})_{xx},
\end{equation}

\begin{equation}\label{s5}
\frac{u_{2}^{n+1} - \widetilde{u_{2}}^{n}}{\triangle t} + F \, \frac{\partial q_{2}^{n}}{\partial u_{1}}\,\frac{u_{1}^{n+1} -  u_{1}^{n}}{\triangle t} + F \, \frac{\partial q_{2}^{n}}{\partial u_{2}}\,\frac{u_{2}^{n+1} -  u_{2}^{n}}{\triangle t}    = D (u_{2}^{n+1})_{xx}.
\end{equation}
The iterative methods in (\ref{s4}) and (\ref{s5})  can be reformulated as
\[
 \frac{\mathbf{u}^{n+1} - \widetilde{\mathbf{u}}^{n}}{ \triangle  t} + F\,  \mathbf{A}^{n}  \frac{{\mathbf{u}}^{n+1} -   \mathbf{u}^{n}}{ \triangle  t}= D \mathbf{u}_{xx}^{n+1},
 \]
 where

\begin{displaymath}
\mathbf{u}^n(x)=
\left(\begin{array}{cc}
u_{1}^n(x)  \\
u_{2}^n(x)\\
\end{array}\right)
\quad  \text{and} \quad \quad \mathbf{A}^{n}=
\left(\begin{array}{cc}
 \frac{\partial q_{1}^{n}}{\partial u_{1}}   &   \frac{\partial q_{1}^{n}}{\partial u_{2}}   \\\\

   \frac{\partial q_{2}^{n}}{\partial u_{1}}  &    \frac{\partial q_{2}^{n}}{\partial u_{2}}   \\
\end{array}\right).
\end{displaymath}

Note that $\frac{\partial q_{1}}{\partial u_{1}}$,$\frac{\partial q_{1}}{\partial u_{2}}$,    $\frac{\partial q_{2}}{\partial u_{1}}$ and  $\frac{\partial q_{2}}{\partial u_{2}}$, are evaluated at the previous time step ($t=t_n$). In order to improve the approximation of (\ref{s4}) and (\ref{s5}) we use the following iteration
\begin{equation}\label{sI4}
\frac{u_{1,l}^{n+1} - \widetilde{u_{1}}^{n}}{\triangle t} + F \, \frac{\partial q_{1,l-1}^{n}}{\partial u_{1}}\,\frac{u_{1,l}^{n+1} -  u_{1}^{n}}{\triangle t }  + F \, \frac{\partial q_{1,l-1}^{n}}{\partial u_{2}}\,\frac{u_{2,l}^{n+1} -  u_{2}^{n}}{\triangle t}    = D (u_{1,l}^{n+1})_{xx},
\end{equation}

\begin{equation}\label{sI5}
\frac{u_{2,l}^{n+1} - \widetilde{u_{2}}^{n}}{\triangle t} + F \, \frac{\partial q_{2,l-1}^{n}}{\partial u_{1}}\,\frac{u_{1,l}^{n+1} -  u_{1}^{n}}{\triangle t} + F \, \frac{\partial q_{2,l-1}^{n}}{\partial u_{2}}\,\frac{u_{2,l}^{n+1} -  u_{2}^{n}}{\triangle t}    = D (u_{2,l}^{n+1})_{xx},
\end{equation}
where $$\frac{\partial q_{1,l}^{n}}{\partial u_{1}}= \frac{ q_{1}(u^{n+1}_{1,l},\cdot)-q_{1}(u_1^{n},\cdot)}{ u_{1,l}^{n+1}-u^n_{1}}.$$
To keep the mass preserved in the   scheme,  we follow the idea of adding perturbation, see\cite{J1}. Define two perturbations of $x_{f}$  by
\[
x^{+}_{f}=  x_{i} -  v \, \triangle t + \eta   \triangle t\, \triangle x,
 \]
 \[
 x^{-}_{f}=  x_{i} -  v \, \triangle t - \eta   \triangle t\, \triangle x,
\]
where the constant $0 < \eta<1  $ depends on $D$, $u$, $ \triangle  t$ and $\triangle x$.  After computing  the values
 $ u_{+}^{n}= u^{n} (x^{+}_{f}) $ and   $u_{-}^{n}=u^{n} (x^{-}_{f}) $ we can compare the amount of  injected concentration for each of components (plus initial concentration if it is not  zero) with the approximated solution  until  level  $n$. If the approximated  mass accumulated up to time level $n$ be less than injected mass, set:
 \[
 u^{n}= \max (u_{+}^{n},  u_{-}^{n}),
\]
otherwise,
\[
 u^{n}= \min (u_{+}^{n},  u_{-}^{n}).
\]

\begin{rem}
One can easily derive the weak formulation and semi-discretized  system and do simulation based on that.
\begin{equation} \label{eq1}
\begin{split}
 \left \langle \mathbf  u^{n} - \tilde{\mathbf u}^{n-1} , \boldsymbol{\phi}\right\rangle +  \left\langle \nabla_{u}   \mathbf q^{n} \, (\mathbf{u^{n}} -\mathbf{u^{n-1}}), \boldsymbol{\phi} \right\rangle +  \triangle t\, \left\langle \nabla\mathbf{u}^{n}, \nabla\boldsymbol{\phi} \right\rangle + \triangle t\, \left\langle \nabla\mathbf{u}^{n}, {\phi} \right\rangle & =  \\
     \triangle t \, \left\langle \mathbf{g}(t^{n},\mathbf{u}^{n-1}), \boldsymbol{\phi}\right\rangle. &
\end{split}
\end{equation}

where
\begin{equation}
\mathbf u^n(x):=\begin{pmatrix}
u_1^n(x)\\
u_2^n(x)
\end{pmatrix}.
\end{equation}
\end{rem}



\subsection{Ideal model}\label{Im}
In the ideal model,  we  assumes that  axial dispersion is negligible i.e., $D=0$  which means that the column has an infinite efficiency and the thermodynamic equilibrium is achieved instantaneously.

\subsubsection{Numerical approach}

We can use MMOCAA explained in the previous section for the ideal case, i.e, $D=0$. Here,  we present a different approach that can be used for one ideal  component ($m=1$), i.e.,
\[
\frac{\partial{u} }{\partial  t } +    F \frac{\partial  }{\partial  t } ( \frac{au}{1+bu })+v\frac{\partial{u}}{\partial  x }  =0.
\]

  Consider the change of variable
\begin{equation}\label{s2}
w= u+ F \frac{a u}{1+bu}.
 \end{equation}
 The idea is  to obtain the approximation of $w$   at point  $(x_{i} , t^{n+1})$. Then  $u$  can  be recovered   as function of  $w$ by the following equation
 \[
u= \frac{-(F a+1-b\, w) + \sqrt{(  Fa+1-b\, w)^2+ 4 bw}}{2b}.
 \]
 By Taylor's expansion we  have
 \begin{equation}\label{fc2}
 w(x_{i}, t^{n+1} )=  w(x_{i}, t^{n} ) +  \triangle t\,  w_{t} (x_{i}, t^{n} )   +  \frac{ (\triangle t)^2 }{2} w_{tt} (x_{i}, t^{n} )  + O (\triangle t)^{3}.
 \end{equation}
Next we obtain approximation for  $w_{t} (x_{i}, t^{n} ) $ and $w_{tt} (x_{i}, t^{n} ). $  To do so, equation (\ref{s2})  implies that:
\begin{equation}\label{c2}
  w_{t}=  (1 +  \frac{ F\, a }{(1+bu)^2} )\, u_{t}.
\end{equation}
By  taking derivative with respect to $t$ from
 \begin{equation}\label{MEQ}
 w_{t} =  - v \, u_{x}
    \end{equation}
  and under some regularity assumption we obtain:
  \[
    w_{tt}  =  (-v\,  u_{x})_{t}=  -v\,(u_{t})_{x}.
 \]
From (\ref{c2}) one has
 \[
 u_{t}= \frac{w_{t}}  {1+   \frac{Fa}{(1+bu)^{2}}}=
 - \frac{ v \, u_{x}}  {1+    \frac{ Fa}{(1+bu)^{2}}}.
 \]
Next  we have
\[
(u_{t})_{x}= -v  (\frac{u_{x}}  {1+  \frac{ F a}{(1+bu)^{2}}})_{x}=  -\frac{v\,  u_{xx}}  {1+ \frac{F a}{(1+bu)^{2}}} - 2 v \frac{F a b (1+bu) u_{x}^2}{( (1+bu)^{2}+ F a)^{2} }.
\]
The recent relation yields
\[
w_{tt} =  \frac{v^2\,  u_{xx}}  {1+ \frac{ F a}{(1+bu)^{2}}} + 2 v^{2} \frac{F a b (1+bu) u_{x}^2}{( (1+bu)^{2}  +F a)^{2}}.
\]
We can substitute $w_{t}$  and $ w_{tt}$ in  (\ref{fc2}) to obtain approximations for $ w(x_{i}, t^{n+1})$ as below:
\begin{equation}\label{LW}
w_i^{n+1}=w_i^n-\frac{  v}{2} \frac{\Delta t }{\Delta x}   \delta_{x}    u_{i}^{n} +\frac{v^2}{2}(\frac{\Delta t}{\Delta x})^2\big(\frac{\delta^2_x u_i^n}{1+\frac{F a}{(1+bu_i^n)^2}}+\frac{2 F ab \,(1+b u_{i}^{n}) ( \delta_{x} u_i^n)^2} {((1+b u_{i}^{n} )^2   + Fa )^2   }\big),
\end{equation}
where $ \delta_{x}$ is the first central difference operator, $\delta^{2}_{x}$ is the second central difference operator, and $\Delta t$ and $\Delta x$ are the mesh-spacing in $t$ and $x$, respectively. The $i$ and $n$ are space and time indices, and $u_i^n$ is the grid function such that $u_i^n = u(x_i,t_n)$.

\section{Numerical implication}

In this section our scheme is validated with different tests. For scalar equation there  exist many approaches   with different flux limiters: Koren, Von leer, superbee, Minmod, Mc. For more detail about this methods refer to \cite{J1}
  \begin{exam}\label{11comp}

To obtain accuracy and compare with analytical solution we
consider the linear adsorption  $q = au.$  The analytical solution of this case with linear adsorption with the parameters given in table \ref{tab1} is derived in \cite{Q}.

\begin{table}[h!]
	 	\caption{Simulation parameters for the  linear case study }
	 	\begin{center}
	 		
	 		\begin{tabular}{|c|c|c|c|}
	 			\hline	
	 			 parameters             &   Symbols            &    Values    &  unite       \\ \hline
	 			 Column length         &  $L$                  &   1          & cm     \\
	 			Porosity               &   $ F$                & 1.5         & -     \\
	 			Interstitial velocity  &  $v$                  &   1      & cm/min       \\
	 			Henry's constant       & $a $                  &  1           & -         \\
	 			 constant in adsorption & $ b$                 & 0             & L/mol         \\
	 			Initial concentration   & $ u_{0} $             &  0              & mol/L          \\
	 			Feed  concentration     &  $ u_{\rm{inj}}$       &  1             & mol/L      \\
                Injection time          &  $ t_{\rm{inj}} $         &  3          & min     \\
	 			Simulation time        &  $ t_{\rm{max}}$              &  7         & min    \\
	 			\hline	
	 		\end{tabular}
	 		\label{tab1}
	 	\end{center}
	 \end{table}

 Figure \ref{figanalnumer} shows both analytical solution and  approximated   solution  with  the   numbers  of spatial steps   $n_x=100$    and of  temporal steps $n_t=400$.

\begin{figure}[h!]
{\includegraphics[width=1.25 \columnwidth]{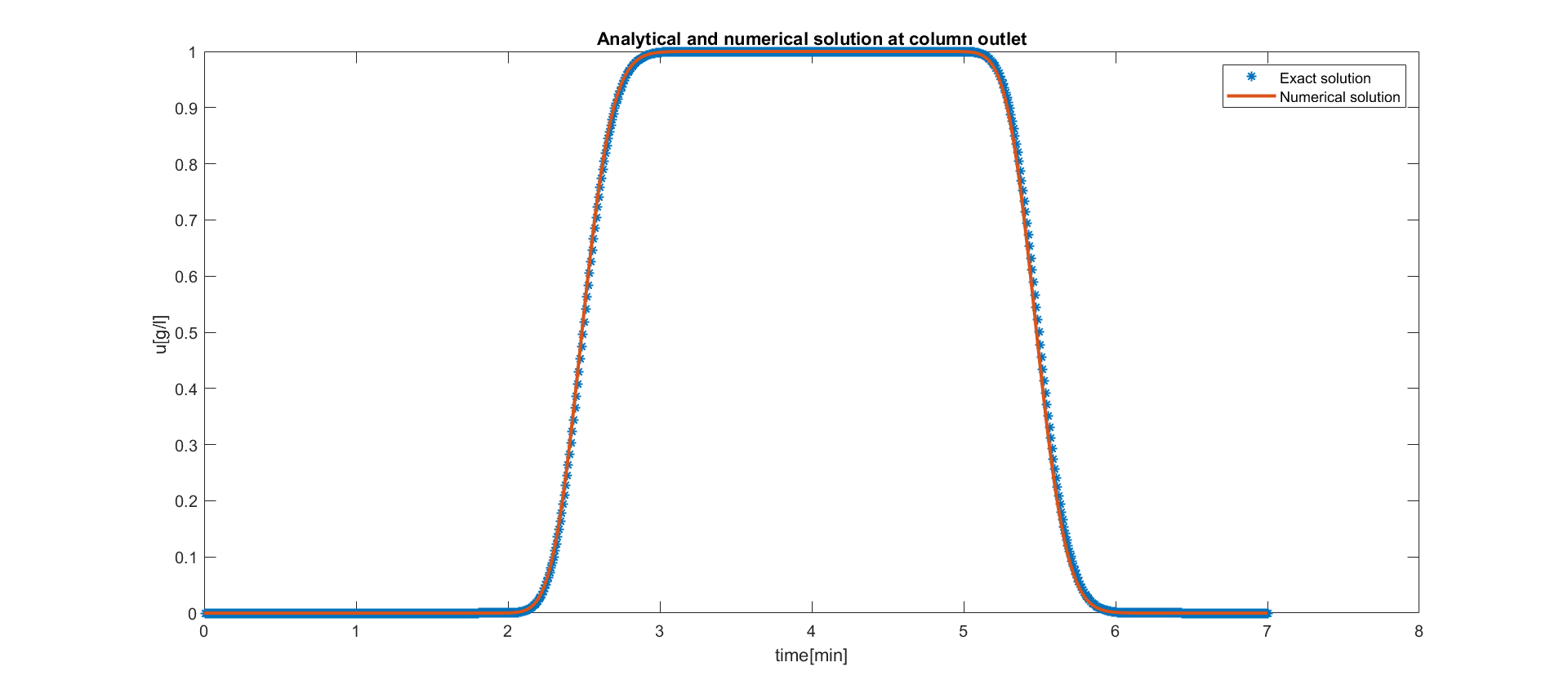}}
\centering
  \caption{Comparison between analytical(continuous red  line)  and approximated solutions at outlet (denoted by  $ \ast$).}
  \label{figanalnumer}
\end{figure}

 Table \ref{tab2} gives a comparison of $L^1$-error and CPU   time
of our method with  discontinuous Galerkin finite element method (DG-FE) with   linear basis functions in \cite{J1, J2} and     with high order basis function of order 8  from \cite{M1}.

\begin{center}
\begin{table}
\scalebox{1}{
\begin{tabular}{ |p{3cm}||p{2cm}|p{2cm}||p{3cm}| }
 \hline
 Different methods    & DOFs  &   $L^1$ error             & CPU time(s)\\\hline
 DG-FM(ord=1) &  16,000&   $0.6 \times 10^{-6}$           &   8827\\\hline
 DG-FM(ord=8) &   90 &    $0.6 \times 10^{-6}$            &0.7\\\hline
 MMOCAA               &  100 &    $0.15\times 10^{-1}$   & 0.11\\\hline
\end{tabular}
}
\caption{The comparison of the method used in \cite{M1} with MMOCAA}\label{tab2}
\end{table}
\end{center}

The $L^1$-norm of  error  and CPU time are presented in  Table \ref{exttable22}.
\begin{center}
\begin{table}
\scalebox{0.9}{
\begin{tabular}{ |p{2cm}|p{2cm}|p{2cm}|p{2.5cm}|  }
 \hline
$n_x$ & $n_t$ & $L^1$ error  & CPU time(s) \\\hline
$50$&$200$&$0.3\times 10^{-1}$&$0.0703$\\\hline
 $100$&$400$&$0.15\times 10^{-1}$&$0.11$\\\hline
$200$&$800$&$0.11\times 10^{-1}$&$0.42$\\\hline
$400$&$1600$&$0.7\times 10^{-2}$&$2.96$\\\hline
$800$&$3200$&$0.5\times 10^{-2}$&$17.64$\\\hline
\end{tabular}
}
\caption{$L^1$-norm and cpu time  for  different $n_x$ and $n_t$.}\label{exttable22}
\end{table}
\end{center}

 \end{exam}


 \begin{exam}\label{1comp}
Here we consider the one component model with nonlinear isotherm given as
\[
q(u)= \frac{u}{1+u}.
\]
The injection time is $0.2 $ and a rectangular pulse of hight  $1\,  g/l$ is injected at inlet.  The length of column
is $1 \textrm{cm,}$ the velocity $v=1\,\textrm{cm/min}$, $ \epsilon =.5.$ and $ N_t=250.$
Figure  \ref{fig1} shows the numerical simulation at outlet, compare with  \cite{J2}.

\begin{figure}[h!]
{\includegraphics[width=.8 \columnwidth]{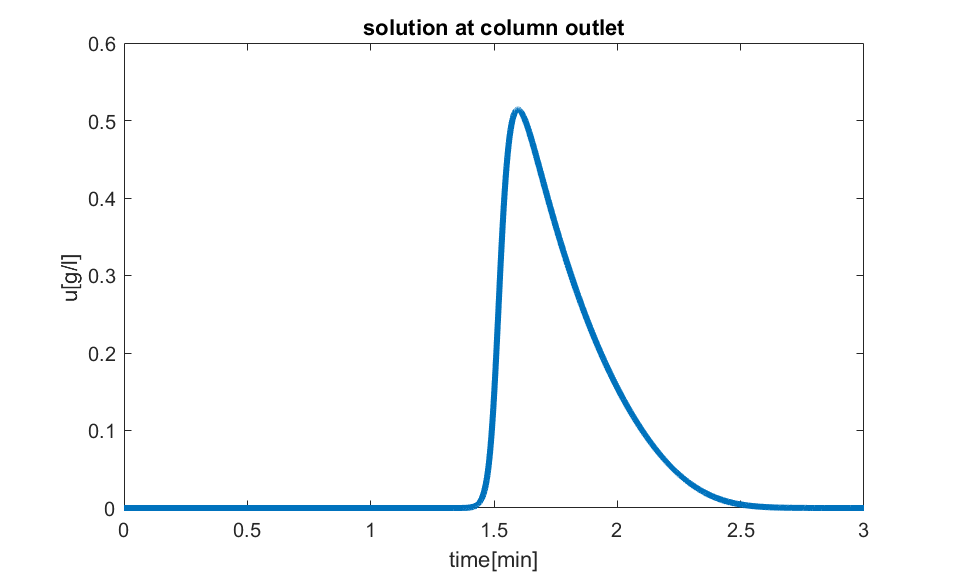}}
\centering
  \caption{Profile of solution $u$ at outlet $x=1$.  }
  \label{fig1}
\end{figure}
We can calculate the mass injected at the  inlet during simulation time. Next, we compute the value of mass passing throughout   each points for the time of simulation.   Figure  \ref{fig2} indicates that the mass is preserved.

\begin{figure}[htb!]
{\includegraphics[width=.8 \columnwidth]{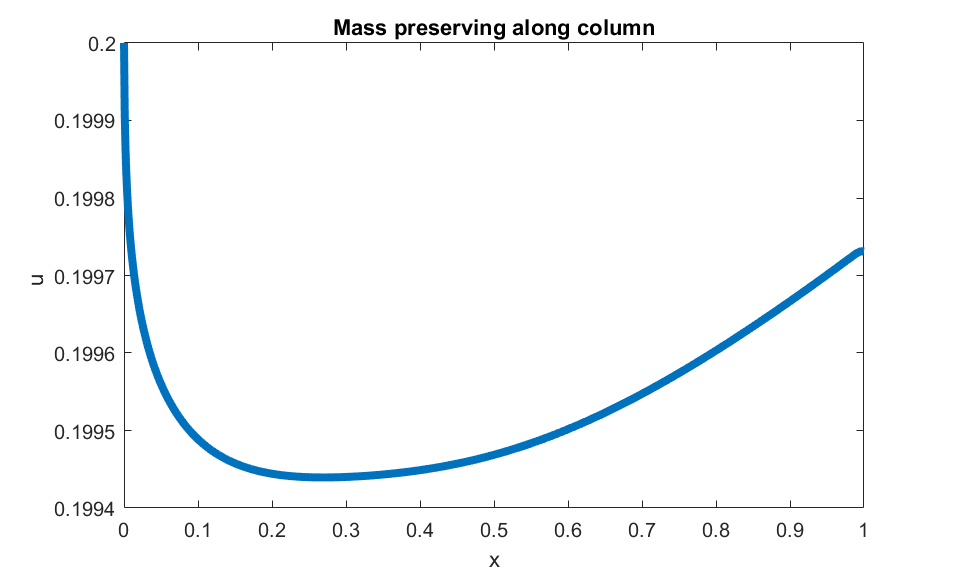}}
\centering
  \caption{The values of approximated solution at different points in column.}
  \label{fig2}
\end{figure}

 Because there is no  analytical solution for this equations  as reference solution,  we consider $n_x=3000 $  grid points and   $n_t =20000$.
 compare the result with the one in \cite{J1}. The $L^1$ error and CPU time are recorded in  Table \ref{tab22}. We compare the results for the case of $n_x= 50$ grid points.
\begin{center}
\begin{table}
\scalebox{0.9}{
\begin{tabular}{ |p{5cm}||p{2cm}|p{2cm}| }
 \hline
 Different methods & $L^1$ error & CPU time(s)\\\hline
First order  &$0.1146$ &0.43\\\hline
Korren   &0.0497 &0.56\\\hline
 Van Leer  & 0.0586 &0.56 \\\hline
 Superbee & 0.0582&0.88 \\\hline
 Minmod     &0.0645&1.45\\\hline
 MC &0.580   &0.62\\\hline
 Our approximation  &  $0.0014 $&$1.5$\\\hline
\end{tabular}
}
\caption{The comparison of the method used in \cite{J1} with  our approximate solution  for one component nonlinear problem}\label{tab22}
\end{table}
\end{center}

\end{exam}


 \begin{exam}\label{2comps}
In this example, we compare our simulation with the test given in  \cite{M1}, section 4.22.  The parameters are chosen from Table \ref{cirdifpq}   with $N_t$= 5000. Here the number of components is two; $m=2$, however, there is  no limitation to simulate   with even larger numbers of theoretical plates. Figure \ref{fig:E} depicts numerical approximation for two components at outlet $x=1$. See Table 6 for $L^1$ error and cpu  time.

 \begin{table}[h]
	 	\caption{Simulation parameters for the nonlinear case study }
	 	\begin{center}
	 		
	 		\begin{tabular}{|c|c|c|c|}
	 			\hline	
	 			 Parameters             &   Symbols            &    Values    &  unite       \\ \hline
	 			 Column length         &  $L$                 &   1           & m     \\
	 			Porosity               &   $ \epsilon $       &  0.4          & -     \\
	 			Interstitial velocity  &  $v$                 &  0.1          & m/s        \\
	 			Henry's constant       & $a_1,  a_2$          &  0.5, 1         & -         \\
	 			 Constant in adsorption & $ b_1, b_2$          & .05, 0.1    & L/mol         \\
	 			Initial concentration   & $ u_{1,0}, u_{2,0}$    &  0, 0         & mol/L          \\
	 			Feed  concentration     &  $ c_{f_1}, c_{f_2}$   &  10, 10         & mol/L      \\
	 			
	 			\hline	
	 		\end{tabular}
	 		\label{cirdifpq}
	 	\end{center}
	 \end{table}

\begin{figure}[h!]
{\includegraphics[width=1.2 \columnwidth]{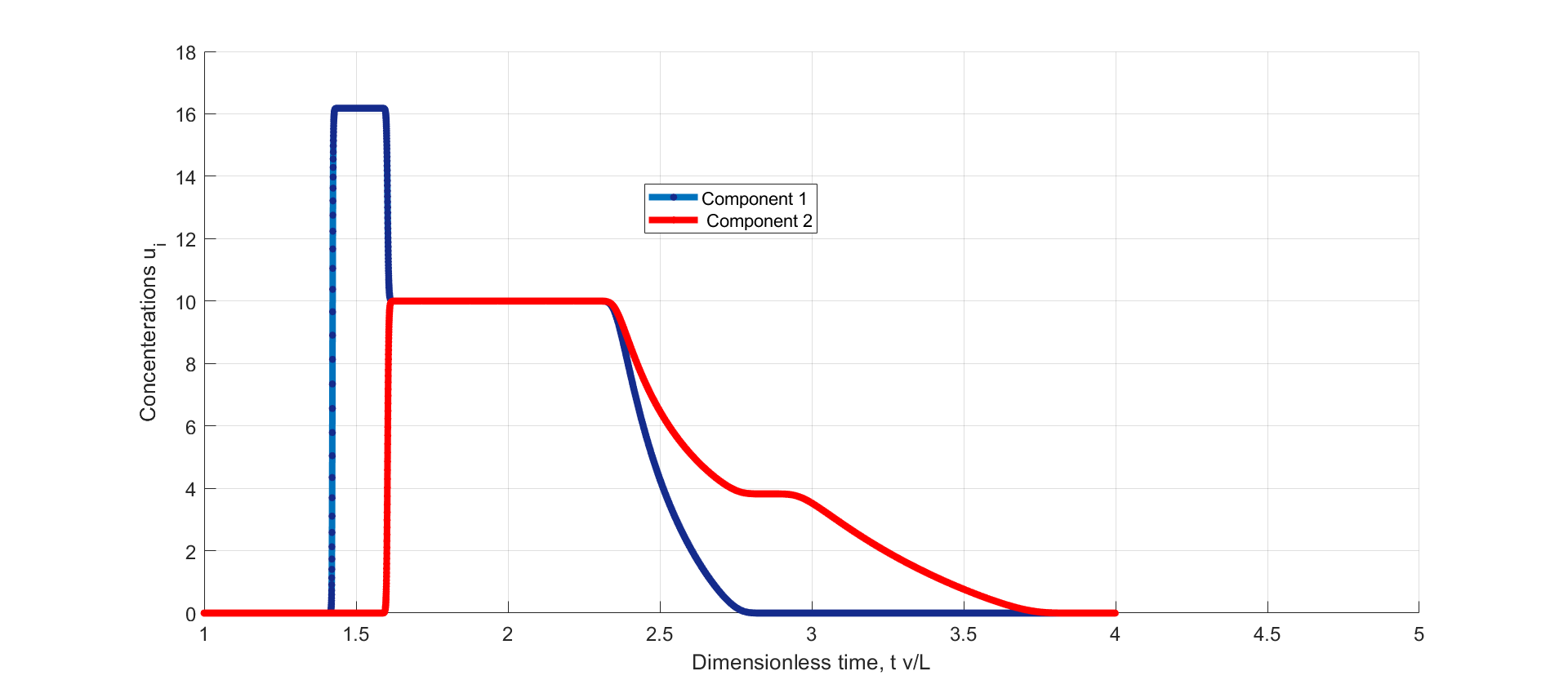}}
\centering
  \caption{Approximate solution for two components at outlet during simulation time.}
  \label{fig:E}
\end{figure}

\begin{center}
\begin{table}
\scalebox{0.9}{
\begin{tabular}{ |p{5cm}||p{2cm}|p{2cm}| }
 \hline
 Different methods & $L^1$ error & CPU time(s)\\\hline
 DG-FM(ord=8) &$0.009$  &4.6\\\hline
 Our approximation         & $0.05  $& 2.3\\\hline
\end{tabular}
}
\caption{The comparison  with   DG-FM in  \cite{M1}}\label{tab222}
\end{table}
\end{center}

\end{exam}

	 \providecommand{\bysame}{\leavevmode\hbox to3em{\hrulefill}\thinspace}
	 \providecommand{\MR}{\relax\ifhmode\unskip\space\fi MR }
	 \providecommand{\MRhref}[2]{%
	 	\href{http://www.ams.org/mathscinet-getitem?mr=#1}{#2}
	 }
	 \providecommand{\href}[2]{#2}

\end{document}